\documentclass[11pt]{amsart}
\usepackage{amsmath,amsthm,amscd,amssymb, color}
\usepackage{latexsym}
\usepackage[colorlinks, citecolor = blue]{hyperref}
\usepackage[capitalize, nameinlink]{cleveref}
\usepackage[alphabetic]{amsrefs}

\usepackage{enumerate}

\usepackage[margin=1.32in]{geometry}

\usepackage{pgfplots}

\newcommand{\Z}{\mathbb Z}

\newcommand{\eps}{\varepsilon}

\newcommand{\avg}{\mathrm{avg}}

\newcommand{\bmx}{\left( \begin{matrix}}
\newcommand{\emx}{\end{matrix} \right)}

\newcommand{\ssk}[2]{ { \Big \{ \begin{matrix} #1 \\ #2 \end{matrix} \Big \} } }

\renewcommand{\mod}{\bmod}

\newtheorem{lem}{Lemma}
\numberwithin{lem}{section}

\newtheorem{conj}[lem]{Conjecture}

\numberwithin{rem}{section}
\newtheorem{question}[lem]{Question}

\numberwithin{equation}{section}

\title{Refined Goldbach conjectures with primes in progressions}

\author{Kimball Martin}
\address{Department of Mathematics, University of Oklahoma, Norman, OK 73019}
\date{\today}

\begin{document}

\maketitle

\begin{abstract} We formulate some refinements of Goldbach's conjectures based
on heuristic arguments and numerical data.  For instance, any even number greater than
4 is conjectured to be a sum of two primes with one prime being 3 mod 4.  In general,
for fixed $m$ and $a, b$ coprime to $m$,
any positive even $n \equiv a + b \bmod m$ outside of a finite exceptional set is expected
to be a sum of two primes $p$ and $q$ with $p \equiv a \bmod m$, $q \equiv b \bmod m$.  
We make conjectures about the growth of these exceptional sets.
\end{abstract}

\section{Introduction}

Let $p$ and $q$ denote prime numbers.  The binary (or strong) Goldbach
conjecture asserts that any even $n > 2$ is of the form $p+q$.  There is both strong heuristic evidence that this is true for sufficiently large $n$ and enormous 
numerical evidence that this is true for $n > 2$.  The same heuristic evidence,
together with equidistibution of primes congruence classes mod $m$, suggests the following:

\begin{conj} \label{conj1} Fix $a, b, m \in \Z$ with $\gcd(a,m) = \gcd(b,m) = 1$.
For sufficiently large even $n \equiv a + b \mod m$, we can write $n=p+q$ for some primes $p \equiv a \mod m$,
$q \equiv b \mod m$.
\end{conj}

While we do not know an explicit statement of this conjecture in the literature,
we do not claim any originality in its formulation.

Denote by $E_{a,b,m}$ the set of positive even $n \equiv a + b \mod m$ which are not
of the form asserted in the conjecture.  This is called the exceptional set for $(a, b, m)$,
and the conjecture asserts $E_{a,b,m}$ is finite.  Note for $a=b=1$ and $m=2$, this is a
(still unknown) weak form of the binary Goldbach conjecture.
Specifically, the binary Goldbach conjecture is equivalent to the statement that
$E_{1,1,2} = \{ 2, 4 \}$.

In this note, we present some heuristic and
numerical investigations on the behavior of these 
exceptional sets.  This leads to both explicit forms of Goldbach's conjecture with 
primes in arithmetic progressions and conjectures about the 
growth of $E_{a,b,m}$.   

First we state a few explicit Goldbach-type conjectures:

\begin{conj} \label{conj:mod4}
Any positive even $n$ as below
is of the form $p+q$ with $p$ and $q$ satisfying the following congruence conditions:

\begin{enumerate}[(i)]
\item any even $n > 4$,  where $p \equiv 3 \mod 4$;

\item any $n \equiv 0 \mod 4$ except $n=4$, where $p \equiv 1 \mod 4$, $q \equiv 3 \mod 4$;

\item any $n \equiv 2 \mod 4$ except $n=2$, where $p \equiv q \equiv 3 \mod 4$;

\item any $n \equiv 2 \mod 4$ except $n=2, 6, 14, 18, 62$, where $p \equiv q \equiv 1 \mod 4$.
\end{enumerate}
\end{conj}

We note that (ii) is already implied by Goldbach's conjecture, whereas (iii) and (iv) are not.
Then (ii) and (iii) imply (i), which may be viewed as a refinement of the binary Goldbach
conjecture.  The notion that the exceptional set should be smaller in case (iii) rather than (iv)
makes sense in light of prime number races, specifically that there are more small primes
which are 3 mod 4 than 1 mod 4.

For moduli $m \ne 4$, 
we just list some sample conjectures about when every (or almost every) multiple of
$m$ is of the form $p+q$ with $p$ and $q$ each coming from single progressions mod $m$:

\begin{conj} \label{conj3}
Any positive even $n$ as below 
is of the form $p+q$ with $p$ and $q$ satisfying the following congruence conditions:

\begin{enumerate}[(i)]

\item any even $n \equiv 0 \mod 3$ except $n=6$, where $p \equiv -q \equiv 1 \mod 3$;

\item any even $n \equiv 0 \mod 5$ (resp.\ except n=$10, 20$), 
where $p \equiv -q \equiv 2 \mod 5$ (resp.\ $p \equiv -q \equiv 1 \mod 5$);

\item any even $n \equiv 0 \mod 7$, where $p \equiv -q \equiv 3 \mod 7$;

\item any even $n \equiv 0 \mod 11$, where $p \equiv -q \equiv 3 \mod 11$;

\item any  $n \equiv 0 \mod 8$, where $p \equiv -q \equiv 3 \mod 8$; 

\item any $n \equiv 0 \mod 16$, where $p \equiv -q \equiv 3 \mod 16$; 

\item any $n \equiv 0 \mod 60$, where $p \equiv -q \equiv a \mod 60$ and 
$a$ is any fixed integer coprime to $60$ which is not $\pm 1, \pm 11 \mod 60$.
\end{enumerate}

\end{conj}

Of these, only the first is a direct consequence of the binary Goldbach conjecture, which we list simply for means of comparison.
The others come from calculations that we describe in
 \cref{sec:numerics}.  

Namely, we compute the following.  First, we may as well assume $m$ is even.
Then we call $(a,b)$ or
$(a, b, m)$ admissible if $a, b \in (\Z/m\Z)^\times$.    
We compute the exceptions $n$ in $E_{a,b,m}$ for all
admissible $(a, b, m)$ with $m \le 200$ up to at least $n = 10^7$.
This appears sufficiently large in our cases of consideration
to believe that we find the full exceptional sets for such $(a, b, m)$.

Likely of more interest than numerous explicit such conjectures is a general
understanding of the behaviour of the exceptional sets.
The first question to ask is:

\begin{question} How fast can $E_{a, b, m}$ grow?
\end{question}

There are a few ways to interpret this: we can look at the growth of the size of
$|E_{a, b, m}|$
or the growth of the sizes of the individual exceptions $n \in E_{a, b, m}$,
and either of these can be interpreted in an average sense or in the sense of looking for an absolute or asymptotic bound.  
Let $E_{\max}(m)$ be the maximum of the exceptions in $E_{a,b,m}$ 
(ranging over $a, b \in (\Z/m\Z)^\times$) and let
$L_\avg(m)$ be the average length (size) of the exceptional sets $E_{a, b, m}$
for a fixed $m$.
Then the heuristics we discuss in \cref{sec:heuristics} suggest the following: 

\begin{conj} \label{conj:asy}
 As $m \to \infty$, we have
\begin{enumerate}[(i)]
\item $E_{\max}(m) :=\max \{ n \in E_{a,b,m} : n \in \Z, \, a, b \in (\Z/m\Z)^\times \} = O(m^2 (\log m)^2)$; and

\item $L_\avg(m) := \frac{1}{\phi(m)^2} \sum_{a, b \in (\Z/m\Z)^\times} |E_{a,b,m}| = 
O(m^\eps)$ for any $\eps > 0$.
\end{enumerate}
\end{conj}

Admittedly our heuristics are rather simplistic (they are too simplistic to
suggest precise asymptotics), but since the numerical data
we present in \cref{sec:numerics} is in strong agreement with these heuristics,
it seems reasonable to believe the above conjecture.  In addition,
our numerical data suggests that $E_{\max}(m)$ grows roughly like a quadratic
function of $m$, which suggests that the growth bound in (i) is not too much
of an overestimate (see \cref{fig2}).

Our data also suggests that as $m$ grows, while the proportion of admissible $(a, b)$
with $E_{a, b, m} = \emptyset$ may decrease on average, there are still many
$(a, b)$ with no exceptional sets (see \cref{tab2}), leading to:
\begin{conj} \label{conj:empty}
There infinitely many tuples $(a, b, m)$ such that
$E_{a,b,m} = \emptyset$.
\end{conj}
We are less confident in this conjecture as precise heuristics about this seem a bit more delicate
(we do not attempt them here),
but it at least seems plausible in connection with \cref{conj:asy}.  Namely, for fixed $m$
the expected length of an exceptional set may be something roughly logarithmic in $m$,
but we have $\phi(m)^2$ such exceptional sets, so there will be a good chance some
of them are empty as long as the variance of $|E_{a,b,m}|$ is not too small.

\medskip
We note that a number of authors (e.g., \cite{lavrik}, \cite{bauer-wang}, \cite{bauer}) have studied the behaviour of
$E_{a,b,m}$ analytically.  However, present analytic methods seem to still be far from
showing finiteness of $E_{a, b, m}$, let alone attacking the finer questions we
explore here.

\medskip
Lastly, we remark that one 
can similarly look at versions of the ternary Goldbach conjecture with primes
in progressions.  However, an answer to the binary case will also give results
about the ternary case, in the same way that the strong Goldbach conjecture
implies the weak Goldbach conjecture (the latter of which is now a theorem \cite{helfgott}).
We simply state one ternary analogue of \cref{conj:mod4}(i):

\begin{conj} Any odd integer $n > 5$ is of the form $n=p+q+r$ for primes $p, q, r$
with $p \equiv q \equiv 2 \mod 3$.
\end{conj}

To see this, our calculations suggest $E_{5,5,6} = \{ 4 \}$ (see \cref{sec:data}), 
which would imply $E_{2,2,3} = \emptyset$,
so in the conjecture we can take $p+q$ to be an arbitrary even number which is 
$1 \mod 3$, and then $r \in \{ 3, 5, 7 \}$.
Unlike the usual ternary Goldbach conjecture, this does not seem to be known at present.  See \cite{li-pan}, \cite{shao}, \cite{shen} for results in this direction.

\subsection*{Acknowledgements}
The author is partially supported by a Simons Foundation Collaboration Grant.

\section{Heuristics} \label{sec:heuristics}

Let $n > 2$ be even, and $g_2(n)$ be the number of ways to write $n=p+q$ for primes
$p$ and $q$.  
In 1922, Hardy and Littlewood conjectured
\[ g_2(n) \sim \mathfrak S(n) \frac{n}{(\log n)^2}, \]
where $\mathfrak S(n)$ is the singular series, which is 0 for odd $n$ and
on average it is is 2 for even $n$.  For a refinement, see \cite{granville}.
For our heuristics, we will naively approximate $g_2(n)$ by the
conjectural average $\frac{2n}{(\log n)^2}$.  This simplification is justified
in our heuristic upper bounds as $\mathfrak S(n) > 1$ for even $n$.

Fix an even modulus $m$.  Let $r$ be the integer part of expected number of admissible
$(a, b) \in ((\Z/m\Z)^\times)^2$
such that $a+b \equiv n$ (averaged over even $0 \le n < m$), 
i.e., $r = [\frac{2\phi(m)^2}{m}]$.
We want to estimate the probability that $n \in E_{a, b, m}$ for some admissible pair
$(a, b) \in ((\Z/m\Z)^\times)^2$.  We will use the following simplistic but reasonable model:
we think of ordered pairs of primes $(p, q)$ solving $p+q=n$ as a collection of $g_2(n)$
independent random events, with the reduction mod $m$ of $(p, q)$ landing in
any of the $r$ admissible classes $(a+m\Z, b + m\Z)$ with equal probability.
(Obviously $(p,q)$ and $(q,p)$ are not independent, but this is not so important for our
heuristics.)
Immediately this suggests \cref{conj1}, but we want to speculate more precisely
on the growth rate of the exceptional sets $E_{a, b, m}$ as $m \to \infty$.

We recall the coupon collector problem.  Say we have $r$ initially empty boxes, 
and at each time $t \in \mathbb N$, a coupon is placed in one box at chosen random.  
Assume at each stage, each box is selected with equal probability $\frac 1r$.
Let $W=W_r$ be the random variable representing the waiting time until all boxes
 have at least 1 coupon.  
 
 Let $X_s$ denote a geometric random variable such that
$P(X_s = k) = (1-s)^{k-1} s$ is the probability of initial success after exactly $k$
trials, where each trial has independent probability of success $s$.
 The problem is to determine the expected value
 $E[W]$.   It is easy to see that $W = X_{r/r} + X_{(r-1)/r} + \cdots X_{1/r}$ (where
 the $X_s$'s are independent), and thus $E[W] = r H_r$, where $H_r = \sum_{j=1}^r \frac 1j$
 is the $r$-th harmonic number.
 
 Thus, in our model, the probability that $n \in E_{a, b, m}$ for some $(a, b)$ is
 simply $P(W_r > g_2(n))$.  One has
 \[ P(W > k) = 1 - \sum_{j=0}^k \frac{r!}{r^j} \ssk {j-1}{r-1} 
 = 1 - \frac{r!}{r^k} \ssk kr = \sum_{j=0}^{r-1} (-1)^{r-j+1} \binom{r}{j}
 \left( \frac j r \right)^k, \]
where $\ssk kr$ denotes the Stirling number of the second kind, i.e., the number of
 ways to partition a set of size $k$ into $r$ nonempty subsets.
  Now we can bound each term  on the right by
 $\binom{r}{[r/2]} (1-1/r)^k$, which will be less than $\frac \eps {k^2 r^4}$ for $r$ large if
 \[ -k \log(1 - \frac 1 r) = k( \frac 1r - \frac 1{r^2} + \cdots ) 
 \gtrsim  r \log 2 + 2 \log k + 4 \log r - \log \eps \gtrsim
  \log \binom{r}{[r/2]} -  \log \frac \eps {k^2 r^4}. \]
 We can make this asymptotic inequality hold by taking $k = Cr^2$, 
 where $C$ depends on $\eps$, and then
 \[ P(W_r > k) < \frac{\eps}{k^2 r^3}. \]
 
 Consider the set $\Sigma_C$ of $\{ (r, k) : r = [\frac{2 \phi(m)^2}m], \,
 k = [ \frac {2n}{(\log n)^2} ] > C r^2 \}$.  Note for $r$ and $k$ of this form,
 the condition $k > Cr^2$ is satisfied when $\frac n{(\log n)^2} > 2 C m^2$,
 and thus when $n > cm^2 (\log m)^2$ for a suitable constant $c$.
Thus, for suitably large $c$, we have a heuristic upper bound on the
 probability that some $E_{a,b,m}$ contains an element $n > cm^2 (\log m)^2$
 of
 \[ \sum_{(r,k) \in \Sigma_C} r P(W_r > k) \le \sum_{(r, k) \in \Sigma_C} \frac{\eps}{k^2 r^2}  < c_0 \eps, \]
 for a uniform constant $c_0$.  The factor of $r$ on the left inside the sum
 comes from accounting for each of the (on average) $r$ classes of pairs $(a, b)$ given $m$, $n$.
 
 This suggests the following bound on the growth of exceptional sets stated
 in \cref{conj:asy}(i):
 \begin{equation}
  E_{\max}(m) = O(m^2 (\log m)^2) \quad \text{as }
 m \to \infty.
 \end{equation}

\medskip
Now let us consider the lengths $|E_{a,b,m}|$ of the exceptional sets.
For fixed $m$, we will model $|E_{a,b,m}|$
as a random variable $L(m)$.
For a fixed $n \equiv a + b \mod m$, the probability (using the model described
above) that $n \in E_{a, b, m}$ is simply $(1 - \frac 1r)^{g_2(n)}$.
Hence the expected size of an exceptional set is
\[ E[L(m)] = \sum_{n \equiv a + b \mod m} \alpha^{g_2(n)}, \quad
\alpha = 1 - \frac 1r. \]
We approximate this with the sum (over all integers $n \ge 2$):
\[ E[L(m)] \approx \frac 1{m} \sum_{n=2}^\infty \alpha^{\frac{2n}{(\log n)^2}}. \]
Then for $0 < \delta < 1$, we have
\[ m E[L(m)] \ll \int_0^\infty \alpha^{2x^\delta} \, dx = \frac 1{2 \delta | \ln \alpha |^{1/\delta} } \int_0^\infty u^{1/\delta - 1} e^{-u} \, du = \frac 1{2 \delta | \ln \alpha |^{1/\delta} } \Gamma ( \frac 1 \delta ). \]
Note $\frac d{dr} {| \ln \alpha |} = \frac d{dr} (\ln (r-1)- \ln r ) = \frac 1{r^2-r}$, so as $r \to
\infty$, we have $\frac 1{|\ln \alpha|} \sim r$.
This gives
\begin{equation} 
E[L(m)] \ll \frac{r^{1/\delta}}{2m} = O(m^\eps), \quad \eps = \frac 1\delta - 1,
\end{equation}
as stated in \cref{conj:asy}(ii).

\section{Numerics} \label{sec:numerics}

Now we present numerical data on the exceptional sets $E_{a, b, m}$ for
(even) $m \le 200$.

\subsection{The method and computational issues}
Our approach, similar to many numerical verifications of Goldbach's conjecture, was roughly
as follows.  To find $\{ n \in E_{a, b, m} : n \le N \}$, we start with two sets of primes
$P = \{ p \equiv a \mod m : p \le M \}$ and $Q = \{ q \equiv b \mod m : q \le N \}$
and determine which $n \le N$ are not of the form $p+q$ for $p \in P$, $q \in Q$,
with $M$ on the order of $10^4$ or $10^5$ depending
on $m$.  Then any potential exceptions below $M$ are guaranteed to actually lie
in $E_{a, b, m}$, and any larger potential exceptions we checked individually by
testing primality (deterministically) of $n-p$ for various $p$.

For each even $m \le 200$, and $a, b \in (\Z/m\Z)^\times$, we checked up to at least
$N = 10^7$ using Sage.  (We checked \cref{conj:mod4} up to $N=10^8$.)
We note that the binary Goldbach conjecture has been numerically
verified for a much, much larger range (up to $4\cdot 10^{18}$ in \cite{oliveira}).  
While one could certainly extend our calculations for larger $N$ (and $m$) with
more efficient implementation and computing resources, our goal here is not to
push the limits of calculation, but rather to generate a reasonable amount of data
to help formulate and support our conjectures.  

That said, there are a couple of obstacles to do a similar amount of verification
for various $E_{a,b,m}$.  First of all, we want to test many triples $(a, b, m)$ which
increases the amount of computation involved.  Second, and much more significant,
when we look for representations $n=p+q$ with $p$ and $q$ in arithmetic progressions,
the minimum value of $p$, say, for which such a representation is possible seems to
increase much faster than without placing congruence conditions on $p$ and $q$.
In other words, to rule out almost all potential exceptions in the first stage
of our algorithm above, for the same $N$ we need to take $M$ larger and larger with $m$.
For instance, when $m=2$ (the usual Goldbach conjecture) one can always take
$M < 10^4$ (i.e., the least prime in a Goldbach partition) to rule out all exceptions for 
$N \le 4 \cdot 10^{18}$ (see \cite{oliveira}), 
however this is not a sufficiently large value of $M$ for many of our calculations.  Already when $N=10^6$ and $M=10^4$,
there are 41 non-exceptions $< 10^6$ that we cannot rule out when $m=50$,  24981 when $m=100$, and 1148651 when $m=148$.

\subsection{Data and observations} \label{sec:data}
For simplicity of exposition, we define
our notation under the hypothesis: \emph{there are no exceptions $n > 10^7$ for $m \le 200$.}
This is believable as the largest exception we find is approximately $10^5$, and
for a given $m$ (with $a, b$ varying) the gaps between one exception and the next
largest exception appear to grow at most quadratically in the number of total exceptions.

Fix $m$.
Let $E_{\max} = E_{\max}(m)$ and $L_{\avg} = L_{\avg}(m)$ be as in the
introduction.  Let $L_{\min} = L_{\min}(m)$ (resp.\  $L_{\max} = L_{\max}(m)$)
be the minimum (resp.\ maximum) of the lengths $|E_{a,b,m}|$ over 
$a, b \in (\Z/m\Z)^\times$.
Let $e_m$ (resp.\ $\tilde e_m$) denote the number of exceptions without
(resp.\ with) multiplicity, i.e., the size of the set (resp.\ multiset) $\bigcup_{(a,b)} E_{a,b,m}$, where $a,b$ run over $(\Z/m\Z)^\times$.

We also consider the above quantities with the additional restriction that
$b \equiv -a \mod m$ so as to treat the special case where $n$ is a multiple of $m$.
In this situation, we denote the analogous quantities with a superscript $0$, e.g.,
$E_{\max}^0$ is the maximal $n \in E_{a, b, m}$ such that $n$ is a multiple of $m$.

We list the first few explicit calculations of exceptional sets (under our hypothesis):

\begin{itemize}
\item $E_{1,1,2} = \{ 2, 4 \}$, which is equivalent to the binary Goldbach conjecture

\item $E_{1,1,4} = \{ 2, 6, 14, 38, 62 \}$, $E_{1, 3, 4} = \{ 4 \}$, and $E_{3, 3, 4} = \{ 2 \}$

\item $E_{1, 1, 6} = \{ 2, 8 \}$, $E_{1, 5, 6} = \{ 6 \}$, and $E_{5, 5, 6} = \{ 4 \}$

\item $E_{1, 1, 8} = \{ 2, 10, 18, 26, 42, 50, 66, 74, 98, 122, 218, 242, 362, 458 \}$,
$E_{1, 3, 8} = \{ 4, 12, 68, 188 \}$, $E_{1, 5, 8} = \{ 6, 14, 38, 62 \}$,
$E_{1, 7, 8} = \{ 8, 16, 32, 56 \}$, $E_{3, 3, 8} = E_{3, 5, 8} = \emptyset$, 
$E_{3, 7, 8} = E_{5, 5, 8} = \{ 2 \}$, $E_{5, 7, 8} = 4$ and $E_{7, 7, 8} = \{ 6, 22, 166 \}$

\item $E_{1, 1, 10} = \{ 2, 12, 32, 152 \}$, $E_{1, 3, 10} = E_{7, 7, 10} = \{ 4 \}$, 
$E_{1, 7, 10} = \{ 8 \}$, $E_{1, 9, 10} = \{ 10, 20 \}$, 
$E_{3, 3, 10} = E_{3, 7, 10} = \emptyset$, $E_{3, 9, 10} = \{ 2, 12 \}$,
$E_{7, 9, 10} = \{ 6, 16 \}$, and $E_{9, 9, 10} = \{ 8, 18, 28, 68 \}$

\end{itemize}

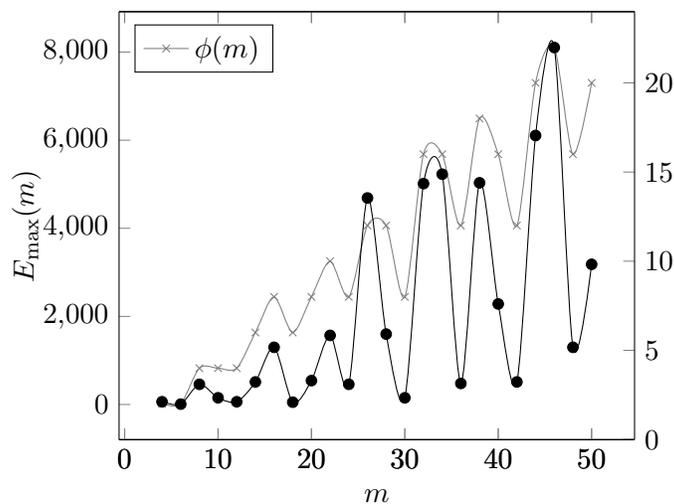
\begin{figure}
\begin{tikzpicture}
    \begin{axis}[xlabel=$m$, axis y line*=right,  legend pos = north west]
    \addplot[smooth,mark=x, gray] plot coordinates {
( 4 , 2 )
( 6 , 2 )
( 8 , 4 )
( 10 , 4 )
( 12 , 4 )
( 14 , 6 )
( 16 , 8 )
( 18 , 6 )
( 20 , 8 )
( 22 , 10 )
( 24 , 8 )
( 26 , 12 )
( 28 , 12 )
( 30 , 8 )
( 32 , 16 )
( 34 , 16 )
( 36 , 12 )
( 38 , 18 )
( 40 , 16 )
( 42 , 12 )
( 44 , 20 )
( 46 , 22 )
( 48 , 16 )
( 50 , 20 )
    };
    \addlegendentry{$\phi(m)$}
    \end{axis}
    \begin{axis}[ylabel = $E_{\max}(m)$, axis y line*=left]
    \addplot[smooth,mark=*]
        plot coordinates {
( 4 , 62 )
( 6 , 8 )
( 8 , 458 )
( 10 , 152 )
( 12 , 62 )
( 14 , 512 )
( 16 , 1298 )
( 18 , 52 )
( 20 , 542 )
( 22 , 1568 )
( 24 , 458 )
( 26 , 4688 )
( 28 , 1598 )
( 30 , 152 )
( 32 , 5014 )
( 34 , 5228 )
( 36 , 478 )
( 38 , 5032 )
( 40 , 2282 )
( 42 , 512 )
( 44 , 6106 )
( 46 , 8104 )
( 48 , 1298 )
( 50 , 3182 )
        };
    \end{axis}
 \end{tikzpicture}
 \caption{Comparing  $E_{\max}(m)$ with $\phi(m)$}
 \label{fig1}
 \end{figure}

We summarize the data from our calculations in \cref{tab1}.
Note that $E_{\max}$ and $E_{\max}^0$, as well as the total number of
exceptions, tend to be relatively larger when
$m$ is a power of 2 or twice a prime.  In these cases $\phi(m)$ is relatively large, i.e.,
we have relatively more admissible pairs $(a, b)$ to consider, so it makes sense that
we pick up more exceptions.  We illustrate this coincidence in the
fluctuations of $E_{\max}(m)$ and $\phi(m)$ in \cref{fig1}.

\begin{table}
\begin{tabular}{r|rrrrrr|rrrrrr}
 $m$ &   $L^0_{\min}$ & $L^0_{\avg}$ &  $L^0_{\max}$ & $E^0_{\max}$ & $e_m^0$ & $\tilde e^0_m$ & $L_{\min}$ & $L_{\avg}$ & $L_{\max}$ & $E_{\max}$& $e_m$ & $\tilde e_m$ \\
\hline
2 & 2 & 2.0 & 2 & 4 & 2 & 2 & 2 & 2.0 & 2 & 4 & 2 & 2 \\
4 & 1 & 1.0 & 1 & 4 & 1 & 2 & 1 & 2.0 & 5 & 62 & 6 & 8 \\
6 & 1 & 1.0 & 1 & 6 & 1 & 2 & 1 & 1.25 & 2 & 8 & 4 & 5 \\
8 & 0 & 2.0 & 4 & 56 & 4 & 8 & 0 & 2.875 & 14 & 458 & 28 & 46 \\
10 & 0 & 1.0 & 2 & 20 & 2 & 4 & 0 & 1.563 & 4 & 152 & 13 & 25 \\
12 & 0 & 0.5 & 1 & 12 & 1 & 2 & 0 & 1.0 & 4 & 62 & 9 & 16 \\
14 & 0 & 1.333 & 3 & 98 & 3 & 8 & 0 & 2.056 & 7 & 512 & 32 & 74 \\
16 & 0 & 2.5 & 4 & 368 & 6 & 20 & 0 & 3.469 & 17 & 1298 & 94 & 222 \\
18 & 0 & 0.333 & 1 & 18 & 1 & 2 & 0 & 0.861 & 2 & 52 & 13 & 31 \\
20 & 0 & 1.5 & 3 & 200 & 4 & 12 & 0 & 1.828 & 8 & 542 & 46 & 117 \\
22 & 0 & 2.2 & 4 & 418 & 7 & 22 & 0 & 2.86 & 9 & 1568 & 102 & 286 \\
24 & 0 & 1.0 & 4 & 192 & 4 & 8 & 0 & 1.344 & 10 & 458 & 39 & 86 \\
26 & 1 & 3.0 & 7 & 754 & 11 & 36 & 0 & 3.313 & 13 & 4688 & 146 & 477 \\
28 & 1 & 2.833 & 7 & 616 & 11 & 34 & 0 & 3.063 & 15 & 1598 & 145 & 441 \\
30 & 0 & 0.25 & 1 & 30 & 1 & 2 & 0 & 0.719 & 3 & 152 & 17 & 46 \\
32 & 0 & 3.875 & 8 & 1184 & 15 & 62 & 0 & 4.969 & 26 & 5014 & 316 & 1272 \\
34 & 1 & 3.5 & 7 & 1088 & 14 & 56 & 0 & 4.355 & 17 & 5228 & 289 & 1115 \\
36 & 0 & 0.833 & 2 & 216 & 4 & 10 & 0 & 1.319 & 5 & 478 & 69 & 190 \\
38 & 2 & 3.889 & 6 & 1558 & 16 & 70 & 0 & 4.864 & 21 & 5032 & 373 & 1576 \\
40 & 1 & 2.375 & 6 & 920 & 9 & 38 & 0 & 2.973 & 19 & 2282 & 225 & 761 \\
42 & 0 & 0.333 & 1 & 42 & 1 & 4 & 0 & 0.896 & 5 & 512 & 39 & 129 \\
44 & 2 & 4.0 & 6 & 3344 & 18 & 80 & 0 & 4.42 & 24 & 6106 & 415 & 1768 \\
46 & 2 & 4.273 & 8 & 1564 & 20 & 94 & 1 & 5.285 & 20 & 8104 & 541 & 2558 \\
48 & 0 & 1.0 & 3 & 288 & 5 & 16 & 0 & 1.629 & 12 & 1298 & 132 & 417 \\
50 & 0 & 2.3 & 4 & 550 & 8 & 46 & 0 & 2.893 & 10 & 3182 & 273 & 1157 \\
52 & 2 & 4.333 & 8 & 3380 & 19 & 104 & 0 & 5.035 & 23 & 8546 & 580 & 2900 \\
54 & 0 & 0.889 & 2 & 216 & 4 & 16 & 0 & 1.503 & 6 & 1096 & 130 & 487 \\
56 & 1 & 3.583 & 6 & 2072 & 12 & 86 & 0 & 4.224 & 24 & 8318 & 491 & 2433 \\
58 & 3 & 6.071 & 13 & 3422 & 29 & 170 & 0 & 6.342 & 31 & 10366 & 870 & 4972 \\
60 & 0 & 0.375 & 2 & 180 & 2 & 6 & 0 & 0.855 & 5 & 542 & 66 & 219 \\
62 & 3 & 6.267 & 10 & 4712 & 28 & 188 & 0 & 6.714 & 28 & 11416 & 975 & 6043 \\
64 & 2 & 6.438 & 12 & 4736 & 32 & 206 & 0 & 7.262 & 35 & 16126 & 1173 & 7436 \\
66 & 0 & 0.9 & 2 & 198 & 3 & 18 & 0 & 1.417 & 6 & 1568 & 141 & 567 \\
68 & 3 & 6.0 & 10 & 5848 & 26 & 192 & 0 & 6.387 & 25 & 13718 & 1044 & 6540 \\
70 & 0 & 1.833 & 3 & 1540 & 9 & 44 & 0 & 2.092 & 9 & 5002 & 280 & 1205 \\
72 & 0 & 1.333 & 4 & 864 & 6 & 32 & 0 & 1.922 & 13 & 2834 & 268 & 1107 \\
74 & 2 & 6.611 & 11 & 6068 & 31 & 238 & 0 & 7.355 & 28 & 16046 & 1370 & 9532 \\
76 & 3 & 6.333 & 11 & 5624 & 30 & 228 & 0 & 6.679 & 29 & 23426 & 1250 & 8656 \\
78 & 0 & 1.0 & 3 & 624 & 5 & 24 & 0 & 1.51 & 7 & 4688 & 205 & 870 \\
80 & 1 & 3.813 & 11 & 3680 & 21 & 122 & 0 & 3.982 & 19 & 6200 & 729 & 4078 \\
82 & 3 & 7.2 & 14 & 8528 & 39 & 288 & 0 & 7.764 & 27 & 25616 & 1651 & 12423 \\
84 & 0 & 0.833 & 4 & 1008 & 5 & 20 & 0 & 1.368 & 11 & 1598 & 202 & 788 \\
86 & 3 & 7.714 & 15 & 12556 & 41 & 324 & 0 & 7.942 & 34 & 26782 & 1805 & 14009 \\
88 & 2 & 5.8 & 13 & 5720 & 29 & 232 & 0 & 6.382 & 32 & 19274 & 1378 & 10211 \\
90 & 0 & 0.333 & 2 & 180 & 2 & 8 & 0 & 1.017 & 5 & 976 & 148 & 586 \\
\end{tabular}
\end{table}

\begin{table}
\begin{tabular}{r|rrrrrr|rrrrrr}
 $m$ &   $L^0_{\min}$ & $L^0_{\avg}$ &  $L^0_{\max}$ & $E^0_{\max}$ & $e_m^0$ & $\tilde e^0_m$ & $L_{\min}$ & $L_{\avg}$ & $L_{\max}$ & $E_{\max}$& $e_m$ & $\tilde e_m$ \\
\hline
92 & 2 & 7.0 & 13 & 7636 & 36 & 308 & 0 & 7.73 & 35 & 21538 & 1891 & 14966 \\
94 & 3 & 8.0 & 16 & 12032 & 39 & 368 & 0 & 8.43 & 38 & 30916 & 2098 & 17838 \\
96 & 0 & 1.375 & 3 & 864 & 6 & 44 & 0 & 2.194 & 17 & 5014 & 437 & 2247 \\
98 & 2 & 5.286 & 10 & 4214 & 25 & 222 & 0 & 5.503 & 28 & 17014 & 1389 & 9708 \\
100 & 0 & 3.8 & 6 & 4300 & 19 & 152 & 0 & 4.366 & 21 & 12134 & 1043 & 6986 \\
102 & 0 & 1.25 & 2 & 918 & 8 & 40 & 0 & 1.848 & 11 & 5228 & 383 & 1892 \\
104 & 3 & 6.708 & 13 & 9256 & 36 & 322 & 0 & 7.236 & 36 & 22886 & 1990 & 16671 \\
106 & 4 & 8.577 & 16 & 11978 & 48 & 446 & 0 & 9.009 & 42 & 39842 & 2661 & 24359 \\
108 & 0 & 2.0 & 4 & 1944 & 13 & 72 & 0 & 2.279 & 13 & 7142 & 548 & 2954 \\
110 & 1 & 2.85 & 8 & 3410 & 15 & 114 & 0 & 3.331 & 14 & 13316 & 822 & 5329 \\
112 & 3 & 5.667 & 12 & 6272 & 29 & 272 & 0 & 6.11 & 34 & 23038 & 1750 & 14077 \\
114 & 0 & 1.278 & 3 & 684 & 6 & 46 & 0 & 2.064 & 12 & 5032 & 482 & 2675 \\
116 & 4 & 8.714 & 15 & 13688 & 44 & 488 & 1 & 8.961 & 42 & 36326 & 2892 & 28102 \\
118 & 4 & 8.793 & 17 & 17228 & 50 & 510 & 0 & 9.446 & 45 & 53614 & 3144 & 31776 \\
120 & 0 & 0.625 & 3 & 360 & 3 & 20 & 0 & 1.3 & 12 & 2282 & 297 & 1331 \\
122 & 4 & 9.4 & 18 & 13298 & 48 & 564 & 0 & 9.601 & 43 & 39818 & 3377 & 34563 \\
124 & 4 & 8.667 & 15 & 12896 & 49 & 520 & 0 & 9.114 & 40 & 42778 & 3265 & 32811 \\
126 & 0 & 1.056 & 4 & 1008 & 7 & 38 & 0 & 1.535 & 11 & 7598 & 389 & 1989 \\
128 & 5 & 9.25 & 15 & 13568 & 47 & 592 & 0 & 10.406 & 43 & 48346 & 3933 & 42622 \\
130 & 1 & 3.458 & 8 & 5590 & 18 & 166 & 0 & 3.955 & 19 & 19406 & 1242 & 9112 \\
132 & 0 & 1.2 & 4 & 1584 & 9 & 48 & 0 & 1.915 & 14 & 6106 & 540 & 3064 \\
134 & 2 & 9.455 & 16 & 14204 & 52 & 624 & 0 & 10.166 & 48 & 52894 & 3992 & 44281 \\
136 & 3 & 8.219 & 17 & 16048 & 47 & 526 & 0 & 8.813 & 52 & 42734 & 3413 & 36100 \\
138 & 0 & 1.682 & 3 & 1656 & 8 & 74 & 0 & 2.332 & 11 & 8104 & 715 & 4515 \\
140 & 0 & 2.792 & 6 & 5180 & 17 & 134 & 0 & 3.18 & 17 & 14198 & 1077 & 7327 \\
142 & 5 & 9.657 & 17 & 23146 & 55 & 676 & 0 & 10.306 & 45 & 54526 & 4328 & 50498 \\
144 & 0 & 2.167 & 5 & 2592 & 15 & 104 & 0 & 2.771 & 18 & 19858 & 959 & 6384 \\
146 & 5 & 10.222 & 22 & 15476 & 54 & 736 & 0 & 10.663 & 54 & 54928 & 4684 & 55278 \\
148 & 2 & 9.528 & 17 & 21608 & 55 & 686 & 0 & 10.128 & 43 & 70222 & 4374 & 52501 \\
150 & 0 & 0.85 & 2 & 1800 & 7 & 34 & 0 & 1.295 & 6 & 3182 & 387 & 2072 \\
152 & 4 & 8.972 & 18 & 17176 & 54 & 646 & 0 & 9.674 & 49 & 62554 & 4331 & 50152 \\
154 & 2 & 4.6 & 9 & 10472 & 25 & 276 & 0 & 4.954 & 22 & 29476 & 1983 & 17836 \\
156 & 0 & 1.625 & 4 & 4212 & 9 & 78 & 0 & 2.204 & 16 & 8546 & 778 & 5078 \\
158 & 6 & 10.718 & 21 & 21488 & 66 & 836 & 1 & 11.011 & 47 & 51248 & 5245 & 66991 \\
160 & 1 & 5.5 & 15 & 13760 & 32 & 352 & 0 & 5.83 & 30 & 31942 & 2479 & 23879 \\
162 & 0 & 2.148 & 7 & 1944 & 12 & 116 & 0 & 2.826 & 15 & 9908 & 1157 & 8241 \\
164 & 5 & 10.6 & 20 & 19352 & 65 & 848 & 0 & 10.768 & 52 & 67546 & 5438 & 68912 \\
166 & 6 & 10.268 & 19 & 22576 & 64 & 842 & 1 & 11.067 & 58 & 65776 & 5817 & 74414 \\
168 & 0 & 1.583 & 8 & 2352 & 11 & 76 & 0 & 1.849 & 15 & 8318 & 675 & 4260 \\
170 & 0 & 4.313 & 10 & 10880 & 27 & 276 & 0 & 4.765 & 31 & 42862 & 2162 & 19517 \\
172 & 4 & 11.024 & 18 & 59168 & 65 & 926 & 0 & 11.098 & 47 & 104494 & 5771 & 78306 \\
174 & 0 & 2.286 & 5 & 7482 & 16 & 128 & 0 & 2.768 & 15 & 10366 & 1166 & 8681 \\
176 & 3 & 8.5 & 16 & 18832 & 53 & 680 & 0 & 8.926 & 52 & 62002 & 4683 & 57129 \\
178 & 5 & 10.955 & 18 & 38092 & 71 & 964 & 1 & 11.528 & 50 & 88622 & 6344 & 89275 \\
180 & 0 & 0.875 & 3 & 1080 & 5 & 42 & 0 & 1.481 & 9 & 7562 & 560 & 3412 \\
\end{tabular}
\end{table}

\begin{table}
\begin{tabular}{r|rrrrrr|rrrrrr}
 $m$ &   $L^0_{\min}$ & $L^0_{\avg}$ &  $L^0_{\max}$ & $E^0_{\max}$ & $e_m^0$ & $\tilde e^0_m$ & $L_{\min}$ & $L_{\avg}$ & $L_{\max}$ & $E_{\max}$& $e_m$ & $\tilde e_m$ \\
\hline
182 & 2 & 5.139 & 14 & 14378 & 33 & 370 & 0 & 5.818 & 30 & 33536 & 2932 & 30161 \\
184 & 5 & 10.023 & 18 & 41768 & 65 & 882 & 1 & 10.665 & 57 & 88106 & 6112 & 82588 \\
186 & 0 & 1.967 & 5 & 2232 & 11 & 118 & 0 & 2.869 & 18 & 11784 & 1285 & 10329 \\
188 & 5 & 10.826 & 17 & 34028 & 70 & 996 & 0 & 11.423 & 55 & 82594 & 6778 & 96686 \\
190 & 2 & 4.639 & 9 & 15200 & 27 & 334 & 0 & 5.052 & 27 & 34652 & 2669 & 26189 \\
192 & 0 & 2.469 & 4 & 7296 & 15 & 158 & 0 & 3.167 & 20 & 16126 & 1552 & 12972 \\
194 & 6 & 11.417 & 24 & 33368 & 75 & 1096 & 1 & 12.054 & 56 & 96728 & 7361 & 111093 \\
196 & 2 & 7.119 & 15 & 20776 & 41 & 598 & 0 & 7.783 & 39 & 56738 & 4514 & 54914 \\
198 & 0 & 1.533 & 4 & 4356 & 9 & 92 & 0 & 2.373 & 16 & 13436 & 1152 & 8542 \\
200 & 2 & 5.6 & 11 & 16400 & 33 & 448 & 0 & 6.228 & 32 & 46922 & 3621 & 39862 \\
\end{tabular}
\caption{Data on exceptional sets mod $m$ for $p \equiv a \mod m$, $q \equiv b \mod m$}
\label{tab1}
\end{table}

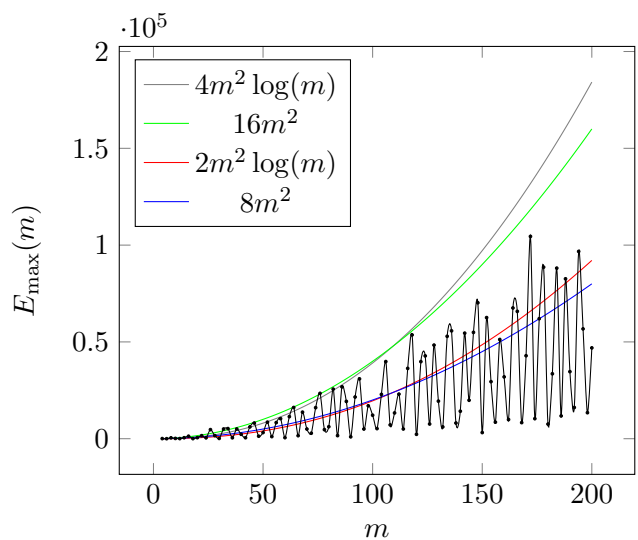
\begin{figure}
\begin{tikzpicture}
    \begin{axis}[xlabel=$m$, ylabel = $E_{\max}(m)$, legend pos = north west]
    \addplot[smooth,mark=., gray]
        plot coordinates {      
( 4 , 11 )
( 6 , 39 )
( 8 , 88 )
( 10 , 160 )
( 12 , 258 )
( 14 , 381 )
( 16 , 532 )
( 18 , 711 )
( 20 , 921 )
( 22 , 1160 )
( 24 , 1431 )
( 26 , 1733 )
( 28 , 2069 )
( 30 , 2437 )
( 32 , 2839 )
( 34 , 3275 )
( 36 , 3745 )
( 38 , 4251 )
( 40 , 4793 )
( 42 , 5370 )
( 44 , 5984 )
( 46 , 6634 )
( 48 , 7322 )
( 50 , 8047 )
( 52 , 8809 )
( 54 , 9610 )
( 56 , 10449 )
( 58 , 11327 )
( 60 , 12244 )
( 62 , 13200 )
( 64 , 14195 )
( 66 , 15230 )
( 68 , 16305 )
( 70 , 17421 )
( 72 , 18576 )
( 74 , 19773 )
( 76 , 21010 )
( 78 , 22289 )
( 80 , 23608 )
( 82 , 24970 )
( 84 , 26372 )
( 86 , 27817 )
( 88 , 29304 )
( 90 , 30833 )
( 92 , 32405 )
( 94 , 34019 )
( 96 , 35676 )
( 98 , 37377 )
( 100 , 39120 )
( 102 , 40906 )
( 104 , 42736 )
( 106 , 44610 )
( 108 , 46527 )
( 110 , 48488 )
( 112 , 50494 )
( 114 , 52543 )
( 116 , 54637 )
( 118 , 56775 )
( 120 , 58958 )
( 122 , 61186 )
( 124 , 63458 )
( 126 , 65776 )
( 128 , 68139 )
( 130 , 70547 )
( 132 , 73000 )
( 134 , 75499 )
( 136 , 78044 )
( 138 , 80634 )
( 140 , 83270 )
( 142 , 85952 )
( 144 , 88680 )
( 146 , 91455 )
( 148 , 94276 )
( 150 , 97143 )
( 152 , 100057 )
( 154 , 103017 )
( 156 , 106024 )
( 158 , 109078 )
( 160 , 112179 )
( 162 , 115327 )
( 164 , 118523 )
( 166 , 121765 )
( 168 , 125055 )
( 170 , 128392 )
( 172 , 131777 )
( 174 , 135209 )
( 176 , 138689 )
( 178 , 142217 )
( 180 , 145793 )
( 182 , 149417 )
( 184 , 153089 )
( 186 , 156809 )
( 188 , 160578 )
( 190 , 164394 )
( 192 , 168260 )
( 194 , 172173 )
( 196 , 176136 )
( 198 , 180147 )
( 200 , 184206 )
};
   \addlegendentry{$4m^2 \log(m)$}

   \addplot[smooth,mark=., green]
        plot coordinates {      
( 4 , 64 )
( 6 , 144 )
( 8 , 256 )
( 10 , 400 )
( 12 , 576 )
( 14 , 784 )
( 16 , 1024 )
( 18 , 1296 )
( 20 , 1600 )
( 22 , 1936 )
( 24 , 2304 )
( 26 , 2704 )
( 28 , 3136 )
( 30 , 3600 )
( 32 , 4096 )
( 34 , 4624 )
( 36 , 5184 )
( 38 , 5776 )
( 40 , 6400 )
( 42 , 7056 )
( 44 , 7744 )
( 46 , 8464 )
( 48 , 9216 )
( 50 , 10000 )
( 52 , 10816 )
( 54 , 11664 )
( 56 , 12544 )
( 58 , 13456 )
( 60 , 14400 )
( 62 , 15376 )
( 64 , 16384 )
( 66 , 17424 )
( 68 , 18496 )
( 70 , 19600 )
( 72 , 20736 )
( 74 , 21904 )
( 76 , 23104 )
( 78 , 24336 )
( 80 , 25600 )
( 82 , 26896 )
( 84 , 28224 )
( 86 , 29584 )
( 88 , 30976 )
( 90 , 32400 )
( 92 , 33856 )
( 94 , 35344 )
( 96 , 36864 )
( 98 , 38416 )
( 100 , 40000 )
( 102 , 41616 )
( 104 , 43264 )
( 106 , 44944 )
( 108 , 46656 )
( 110 , 48400 )
( 112 , 50176 )
( 114 , 51984 )
( 116 , 53824 )
( 118 , 55696 )
( 120 , 57600 )
( 122 , 59536 )
( 124 , 61504 )
( 126 , 63504 )
( 128 , 65536 )
( 130 , 67600 )
( 132 , 69696 )
( 134 , 71824 )
( 136 , 73984 )
( 138 , 76176 )
( 140 , 78400 )
( 142 , 80656 )
( 144 , 82944 )
( 146 , 85264 )
( 148 , 87616 )
( 150 , 90000 )
( 152 , 92416 )
( 154 , 94864 )
( 156 , 97344 )
( 158 , 99856 )
( 160 , 102400 )
( 162 , 104976 )
( 164 , 107584 )
( 166 , 110224 )
( 168 , 112896 )
( 170 , 115600 )
( 172 , 118336 )
( 174 , 121104 )
( 176 , 123904 )
( 178 , 126736 )
( 180 , 129600 )
( 182 , 132496 )
( 184 , 135424 )
( 186 , 138384 )
( 188 , 141376 )
( 190 , 144400 )
( 192 , 147456 )
( 194 , 150544 )
( 196 , 153664 )
( 198 , 156816 )
( 200 , 160000 )
};
   \addlegendentry{$16m^2$}

   \addplot[smooth,mark=., red]
        plot coordinates {   
   ( 4 , 5 )
( 6 , 19 )
( 8 , 44 )
( 10 , 80 )
( 12 , 129 )
( 14 , 190 )
( 16 , 266 )
( 18 , 355 )
( 20 , 460 )
( 22 , 580 )
( 24 , 715 )
( 26 , 866 )
( 28 , 1034 )
( 30 , 1218 )
( 32 , 1419 )
( 34 , 1637 )
( 36 , 1872 )
( 38 , 2125 )
( 40 , 2396 )
( 42 , 2685 )
( 44 , 2992 )
( 46 , 3317 )
( 48 , 3661 )
( 50 , 4023 )
( 52 , 4404 )
( 54 , 4805 )
( 56 , 5224 )
( 58 , 5663 )
( 60 , 6122 )
( 62 , 6600 )
( 64 , 7097 )
( 66 , 7615 )
( 68 , 8152 )
( 70 , 8710 )
( 72 , 9288 )
( 74 , 9886 )
( 76 , 10505 )
( 78 , 11144 )
( 80 , 11804 )
( 82 , 12485 )
( 84 , 13186 )
( 86 , 13908 )
( 88 , 14652 )
( 90 , 15416 )
( 92 , 16202 )
( 94 , 17009 )
( 96 , 17838 )
( 98 , 18688 )
( 100 , 19560 )
( 102 , 20453 )
( 104 , 21368 )
( 106 , 22305 )
( 108 , 23263 )
( 110 , 24244 )
( 112 , 25247 )
( 114 , 26271 )
( 116 , 27318 )
( 118 , 28387 )
( 120 , 29479 )
( 122 , 30593 )
( 124 , 31729 )
( 126 , 32888 )
( 128 , 34069 )
( 130 , 35273 )
( 132 , 36500 )
( 134 , 37749 )
( 136 , 39022 )
( 138 , 40317 )
( 140 , 41635 )
( 142 , 42976 )
( 144 , 44340 )
( 146 , 45727 )
( 148 , 47138 )
( 150 , 48571 )
( 152 , 50028 )
( 154 , 51508 )
( 156 , 53012 )
( 158 , 54539 )
( 160 , 56089 )
( 162 , 57663 )
( 164 , 59261 )
( 166 , 60882 )
( 168 , 62527 )
( 170 , 64196 )
( 172 , 65888 )
( 174 , 67604 )
( 176 , 69344 )
( 178 , 71108 )
( 180 , 72896 )
( 182 , 74708 )
( 184 , 76544 )
( 186 , 78404 )
( 188 , 80289 )
( 190 , 82197 )
( 192 , 84130 )
( 194 , 86086 )
( 196 , 88068 )
( 198 , 90073 )
( 200 , 92103 )
};
   \addlegendentry{$2m^2 \log(m)$}

   \addplot[smooth,mark=., blue]
        plot coordinates {      
( 4 , 32 )
( 6 , 72 )
( 8 , 128 )
( 10 , 200 )
( 12 , 288 )
( 14 , 392 )
( 16 , 512 )
( 18 , 648 )
( 20 , 800 )
( 22 , 968 )
( 24 , 1152 )
( 26 , 1352 )
( 28 , 1568 )
( 30 , 1800 )
( 32 , 2048 )
( 34 , 2312 )
( 36 , 2592 )
( 38 , 2888 )
( 40 , 3200 )
( 42 , 3528 )
( 44 , 3872 )
( 46 , 4232 )
( 48 , 4608 )
( 50 , 5000 )
( 52 , 5408 )
( 54 , 5832 )
( 56 , 6272 )
( 58 , 6728 )
( 60 , 7200 )
( 62 , 7688 )
( 64 , 8192 )
( 66 , 8712 )
( 68 , 9248 )
( 70 , 9800 )
( 72 , 10368 )
( 74 , 10952 )
( 76 , 11552 )
( 78 , 12168 )
( 80 , 12800 )
( 82 , 13448 )
( 84 , 14112 )
( 86 , 14792 )
( 88 , 15488 )
( 90 , 16200 )
( 92 , 16928 )
( 94 , 17672 )
( 96 , 18432 )
( 98 , 19208 )
( 100 , 20000 )
( 102 , 20808 )
( 104 , 21632 )
( 106 , 22472 )
( 108 , 23328 )
( 110 , 24200 )
( 112 , 25088 )
( 114 , 25992 )
( 116 , 26912 )
( 118 , 27848 )
( 120 , 28800 )
( 122 , 29768 )
( 124 , 30752 )
( 126 , 31752 )
( 128 , 32768 )
( 130 , 33800 )
( 132 , 34848 )
( 134 , 35912 )
( 136 , 36992 )
( 138 , 38088 )
( 140 , 39200 )
( 142 , 40328 )
( 144 , 41472 )
( 146 , 42632 )
( 148 , 43808 )
( 150 , 45000 )
( 152 , 46208 )
( 154 , 47432 )
( 156 , 48672 )
( 158 , 49928 )
( 160 , 51200 )
( 162 , 52488 )
( 164 , 53792 )
( 166 , 55112 )
( 168 , 56448 )
( 170 , 57800 )
( 172 , 59168 )
( 174 , 60552 )
( 176 , 61952 )
( 178 , 63368 )
( 180 , 64800 )
( 182 , 66248 )
( 184 , 67712 )
( 186 , 69192 )
( 188 , 70688 )
( 190 , 72200 )
( 192 , 73728 )
( 194 , 75272 )
( 196 , 76832 )
( 198 , 78408 )
( 200 , 80000 )
};
   \addlegendentry{$8m^2$}

    \addplot[smooth,mark=*, mark size = 0.5]
        plot coordinates {
( 4 , 62 )
( 6 , 8 )
( 8 , 458 )
( 10 , 152 )
( 12 , 62 )
( 14 , 512 )
( 16 , 1298 )
( 18 , 52 )
( 20 , 542 )
( 22 , 1568 )
( 24 , 458 )
( 26 , 4688 )
( 28 , 1598 )
( 30 , 152 )
( 32 , 5014 )
( 34 , 5228 )
( 36 , 478 )
( 38 , 5032 )
( 40 , 2282 )
( 42 , 512 )
( 44 , 6106 )
( 46 , 8104 )
( 48 , 1298 )
( 50 , 3182 )
( 52 , 8546 )
( 54 , 1096 )
( 56 , 8318 )
( 58 , 10366 )
( 60 , 542 )
( 62 , 11416 )
( 64 , 16126 )
( 66 , 1568 )
( 68 , 13718 )
( 70 , 5002 )
( 72 , 2834 )
( 74 , 16046 )
( 76 , 23426 )
( 78 , 4688 )
( 80 , 6200 )
( 82 , 25616 )
( 84 , 1598 )
( 86 , 26782 )
( 88 , 19274 )
( 90 , 976 )
( 92 , 21538 )
( 94 , 30916 )
( 96 , 5014 )
( 98 , 17014 )
( 100 , 12134 )
( 102, 5228 )
( 104, 22886)
( 106, 39842)
( 108, 7142)
( 110, 13316)
( 112, 23038)
( 114, 5032)
( 116, 36326)
( 118,53614 )
( 120, 2282)
( 122, 39818)
( 124, 42778)
( 126, 7598)
( 128, 48346)
( 130, 19406)
( 132, 6106)
( 134, 52894)
( 136, 55714)
( 138, 8104)
( 140, 14198)
( 142, 54526)
( 144, 19858)
( 146, 54928)
( 148, 70222)
( 150, 3182)
( 152, 62554)
( 154, 29476)
( 156, 8546)
( 158, 51248)
( 160, 31942)
( 162, 9908)
( 164, 67546)
( 166, 65776)
( 168, 8318)
( 170, 42862)
( 172, 104494)
( 174, 10366)
( 176, 62002)
( 178, 88622)
( 180, 7562)
( 182, 33536)
( 184, 88106)
( 186, 11784)
( 188, 82594)
( 190, 34652)
( 192, 16126)
( 194, 96728)
( 196, 56738)
( 198, 13436)
( 200, 46922)
        };   
   
    \end{axis}
 \end{tikzpicture}
 \caption{Comparing $E_{\max}(m)$ with quadratic and quadratic times logarithmic growth}
 \label{fig2}
 \end{figure}
 
 \begin{table}
 \begin{tabular}{r||rrrrr rrrrr rrrrr}
$m$ & 2 & 4 & 6 & 8 & 10 & 12 & 14 & 16 & 18 & 20 & 22 & 24 & 26 & 28 & 30 \\
\hline
$\#$ & 0 & 0 & 0 & 3 & 3 & 3 & 6 & 6 & 9 & 11 & 2 & 20 & 10 & 16 & 21 \\ 
$\%$ &  0.0 & 0.0 & 0.0 & 18.8 & 18.8 & 18.8 & 16.7 & 9.4 & 25.0 & 17.2 & 2.0 & 31.3 & 6.9 & 11.1 & 32.8 \\
\hline
\hline
$m$ & 32 & 34 & 36 & 38 & 40 & 42 & 44 & 46 & 48 & 50 & 52 & 54 & 56 & 58 & 60 \\
\hline
$\#$ & 12 & 12 & 33 & 11 & 26 & 44 & 8 & 0 & 68 & 32 & 10 & 75 & 18 & 4 & 97 \\
$\%$ &
4.7 & 4.7 & 22.9 & 3.4 & 10.2 & 30.6 & 2.0 & 0.0 & 26.6 & 8.0 & 1.7 & 23.1 & 3.1 & 0.5 & 37.9 \\
\hline
\hline
$m$ & 62 & 64 & 66 & 68 & 70 & 72 & 74 & 76 & 78 & 80 & 82 & 84 & 86 & 88 & 90 \\
\hline
$\#$ & 8 & 4 & 108 & 24 & 60 & 105 & 4 & 14 & 118 & 42 & 6 & 160 & 14 & 2 & 195 \\
$\%$ & 0.9 & 0.4 & 27.0 & 2.3 & 10.4 & 18.2 & 0.3 & 1.1 & 20.5 & 4.1 & 0.4 & 27.8 & 0.8 & 0.1 & 33.9 \\
\hline
\hline
$m$ & 92 & 94 & 96 & 98 & 100 & 102 & 104 & 106 & 108 & 110 & 112 & 114 & 116 & 118 & 120 \\
\hline
$\#$ & 6 & 14 & 163 & 26 & 18 & 147 & 20 & 6 & 171 & 84 & 12 & 173 & 0 & 6 & 326 \\
$\%$ & 0.3 & 0.7 & 15.9 & 1.5 & 1.1 & 14.4 & 0.9 & 0.2 & 13.2 & 5.2 & 0.5 & 13.3 & 0.0 & 0.2 & 31.8 \\
\hline
\hline
$m$ & 122 & 124 & 126 & 128 & 130 & 132 & 134 & 136 & 138 & 140 & 142 & 144 & 146 & 148 & 150 \\
\hline
$\#$ & 6 & 20 & 286 & 4 & 64 & 297 & 4 & 8 & 194 & 107 & 2 & 241 & 6 & 2 & 422 \\
$\%$ & 0.2 & 0.6 & 22.1 & 0.1 & 2.8 & 18.6 & 0.1 & 0.2 & 10.0 & 4.6 & 0.0 & 10.5 & 0.1 & 0.0 & 26.4 \\
\hline
\hline
$m$ & 152 & 154 & 156 & 158 & 160 & 162 & 164 & 166 & 168 & 170 & 172 & 174 & 176 & 178 & 180 \\
\hline
$\#$ & 4 & 71 & 272 & 0 & 95 & 240 & 6 & 0 & 473 & 100 & 4 & 314 & 8 & 0 & 609 \\
$\%$ & 0.1 & 2.0 & 11.8 & 0.0 & 2.3 & 8.2 & 0.1 & 0.0 & 20.5 & 2.4 & 0.1 & 10.0 & 0.1 & 0.0 & 26.4 \\
 \end{tabular}
 \caption{Counting the number $(a,b)$ for which $|E_{a,b,m}| = 0$}
 \label{tab2}
 \end{table}

While our numerics are somewhat limited, they suggest that the growth of 
$E_{\max}(m)$ appears to be strictly slower than that of $m^2 (\log m)^2$ as
stated in \cref{conj:asy}, and the true growth rate appears to be closer
to $O(m^2)$ or $O(m^2 \log m)$---see \cref{fig2} for an overlay
of the graphs of $E_{\max}(m)$ (black dots, with the scale on the left) and $\phi(m)$ (
gray x's, with the scale on the right).

Finally, we observe it often happens that (under our hypothesis), for fixed $m$,
at least one admissible pair $(a, b)$ has no exceptions, i.e., every
$n \equiv a + b \mod m$ is of the form $p+q$ with $p \equiv a \mod m$ and
$q \equiv b \mod m$.  Both the number of such pairs $(a, b)$ and the fraction
of such pairs out of the total number $\phi(m)^2$ of admissible pairs
are tabulated in \cref{tab2}.  It appears that
$z_m = \# \{ (a, b) \in ((\Z/m\Z)^\times)^2 : E_{a, b, m} = \emptyset \}$ tends
to grow at least on average in several of the columns in \cref{tab2} 
(e.g., when $m$ is a multiple of 5 or 6).  This suggests that $z_m$ is unbounded,
and in particular suggests \cref{conj:empty}.


\end{document}